\documentclass[a4 paper,11pt,twoside]{article}
\usepackage[centertags]{amsmath}
\usepackage{amssymb}
\usepackage{amsfonts}
\usepackage{amsthm}
\usepackage{newlfont}
\usepackage{amscd}
\usepackage{tikz}
\usepackage[top=1.7in, bottom=1.7in, left=1.5in, right=1.5in]{geometry}
\title{\large FUZZY EXTENDED FILTERS OF $\emph{MS}$-ALGEBRAS}
\author{}
\date{}
\theoremstyle{plain}
\newtheorem{theorem}{Theorem}[section]
\newtheorem{proposition}[theorem]{Proposition}
\newtheorem{lemma}[theorem]{Lemma}
\newtheorem{corollary}[theorem]{Corollary}
\newtheorem{remark}[theorem]{Remark}
\newtheorem{example}[theorem]{Example}
\theoremstyle{definition}
\newtheorem{definition}[theorem]{Definition}
\begin{document}\maketitle
\begin{center}
\textbf{Ahmed Gaber$^{1}$, M.A.Seoud$^{2}$and Mona Tarek$^{3}$}\vspace{0.15 cm}\\
\small{$^{1,2,3}$Department of Mathematics\\Faculty of Science, Ain Shams University, Cairo, Egypt\\
$^{1}$a.gaber@sci.asu.edu.eg\vspace{0.15cm}, 0000-0002-5620-901X\\  
$^{2}$m.a.seoud@hotmail.com\vspace{0.15cm}\\ 
$^{3}$Mona.Saad@sci.asu.edu.eg\vspace{0.15cm}, 0000-0002-0838-2563\\ }
\end{center}
\vspace{0.1 cm}

\begin{abstract}
  In this article, for an $\emph{MS}$-
  algebra and a fuzzy filter $\chi $, the 
  concept of extended fuzzy filter of  $\chi$ is 
  presented, notated by $
  \Upsilon_{\chi,W}$ with $W \subseteq 
  \mathcal{L}$. The 
  features of $\Upsilon_{\chi,W}$ are 
  investigated. Furthermore, the strong 
  fuzzy filter is introduced, donated by 
  $\Omega_{\chi,W}$.
  Many properties are studied.  
  Characterisation of both $
  \Upsilon_{\chi,W}, \Omega_{\chi,W}$ are clarified    
  by using the notion of dense elements with 
  respect to a fuzzy filter. The homomorphisms 
  of both $\Upsilon_{\chi,W}$ and $\Omega_{\chi,W}$ are
  triggered. 
\end{abstract}
\begin{flushleft}
\end{flushleft}
\begin{flushleft}
\end{flushleft}
\section{\textbf{I}ntroduction}
\hspace*{0.4 cm}
In 1965, the genesis of the fuzzy mathematics 
was made by Lotfi Asker Zadeh'publications 
[15], [16] and [17]. Lotfi Asker Zadeh 
introduced several fuzzy concepts 
as fuzzy sets, fuzzy algorithms, etc. 
Many algebraic structures were 
fuzzified. In 1971, fuzzy subgroups 
and fuzzy subgroupoids were settled down 
by A.Rosenfeld [11].   
The evolution of fuzzy mathematics took a 
great part in lattice theory see[1], [2], 
[3] and [10].\\
\hspace*{0.4 cm}
The variety $\mathbf{MS}$ that contains all $
\emph{MS}$-algebras was first considered 
by T.S. Blyth and J.C. Varlet to provide a  
unified frame work to handle both de 
Morgan algebras and Stone algebras [4], [5], [6] and [7]. The 
 theory of substructures of $\emph{MS}$-algebras has been studied by many
 researchers 
 as in [12], [13] and [14].\\
  \hspace*{0.4 cm}   
In this paper, for an $\emph{MS}$-algebra 
$\mathcal{L}$ and a fuzzy filter $\chi$ 
together with $W\subseteq \mathcal{L}$, 
the fuzzification of the extended filter  
$E_\chi(W)$ is presented and notated by $
  \Upsilon_{\chi,W}$. We verify that $
  \Upsilon_{\chi,W}$ is an extended prime fizzy 
  filter for $\chi$ and properties of $
  \Upsilon_{\chi,W}$ are studied in a 
  great details. The class of all  
extended fuzzy filters is proved to be  
closed under the meet operation. The 
notion of a strong extended filter 
and notated by $\Omega_{\chi,W}$. The 
values of both $\Upsilon_{\chi,W} 
, \Omega_{\chi,W}$ are 
studied in terms of dense elements. The features of $
\Upsilon_{\chi,W}$ 
homomorphism are characterised. We study 
the elements of the kernel and the 
cokernel.   
\section{\textbf{P}reliminaries}
\hspace*{0.4 cm} To make the paper 
consistent, 
we mention some results and definitions 
that will be used throughout the paper.
For $\chi \in \mathcal{F}(\mathcal{L})$.
We call $\chi \subseteq \mathcal{L}$ a 
filter providing that $\chi$ is a 
sublattice and if $ 
 \theta \in \chi$, $ \varepsilon \in \mathcal{L}$ imply 
 that $\theta \vee \varepsilon \in \chi$. 
 A proper filter $\mathcal{P}$ is 
 promoted to be prime if $\theta,
 \varepsilon \in \mathcal{L}$ 
 such that $ \theta \vee \varepsilon \in 
 \mathcal{P}$, then $ \theta \in 
 \mathcal{P}$ or $ \varepsilon \in 
 \mathcal{P}$. Let $\varepsilon\in 
 \mathcal{L}$. The 
 principal filter $[\varepsilon)$ of 
 $\mathcal{L}$ generated by $\varepsilon$ 
 is equivalent to $[\varepsilon)=\{\theta\in \mathcal{L}:
 \theta\geq \varepsilon
 \}$. For every $W\subseteq \mathcal{L}$, 
 the filter $[\mathfrak{B})$ generated by the set $\mathfrak{B}$ is defined to be 
 \begin{center}
$[\mathfrak{B})=\{\theta\in \mathcal{L}:\theta\geq 
b_1\wedge b_2\wedge...\wedge b_n$ for 
$b_1,b_2,...,b_n\in \mathfrak{B} \}$. 
\end{center}
 For a lattice $\mathcal{L}$ endowed with the distributive property, the filters of
 $\mathcal{L}$ ordered by inclusion is a
 lattice which denoted by $\mathcal{F}(\mathcal{L})$. For a bounded distributive lattice $\mathcal{L}$, the filter $[1)=\{1\}$ is the smallest member of $\mathcal{F}(\mathcal{L})$ and $[0)=\mathcal{L}$ is the largest member of $\mathcal{F}(\mathcal{L})$.
\vspace{0.1 cm}
\newline \hspace*{0.4 cm}
If $\theta\mapsto \theta^{\circ}$ is a 
unary operation on a bounded distributive lattice $(\mathcal{L};\vee,\wedge, 0, 1)$
 satisfying,
\begin{center}
$1^{\circ}=0$, $(\ \theta\wedge \xi \ )^{\circ}
=\theta^{\circ}\vee
\xi^{\circ}$ and $\theta \leq 
\theta^{\circ\circ}$ .
\end{center}
Then $\mathcal{L}$ is called an $\emph{MS}$-algebra, see $[6]$.
The elements in an $\emph{MS}$-algebra satisfy the following equalities.
Throughout the paper we use the notation $\mathcal{L}$ for an $\emph{MS}$-algebra.
We refer to $[9]$ for details of the concepts of lattice theory.
\begin{proposition}$[8]$
Let $\theta, \xi \in \mathcal{L}$. The following 
equalities are valid; 
\item{(1)} $(\ \xi\vee \theta\ )^{\circ}=
\xi^{\circ}\wedge \theta^{\circ}$,
\item{(2)} $(\ \xi \vee \theta \ )^{\circ\circ}=\xi^{\circ\circ}\vee \theta^{\circ\circ}$,
\item{(3)} $\theta^{\circ\circ\circ}=\theta^{\circ}$,
\item{(4)} $0^{\circ}=1$.
\end{proposition}
 
    \begin{definition}
    Let $\chi \in\mathcal{F}(\mathcal{L})$ and
    $\mathcal{L} \in \mathbf{MS}$. If $W$ 
    is a 
    nonempty subset of $\mathcal{L}$, define
\begin{center}
    $E_\chi(W)=\{ \theta \in \mathcal{L}\  ;\theta\vee w^{\circ \circ} \in \chi \ for \ every \ w \in W \}$ .
\end{center}
    \end{definition}
    Recall some basic theorems that characterise the 
    concept of $E_\chi(W)$.
    \begin{theorem}
Let $\chi \in \mathcal{F}(\mathcal{L})$, then $E_\chi(W)$ is  a filter of $L$ containing $\chi$. 
\end{theorem}
\begin{proof}
Obviously, $1 \in E_\chi(W)$. Assume that $\theta \in E_\chi(W)$ 
and $\varrho\in \mathcal{L}$ satisfying $\theta \leq \varrho$. We have that $\varrho\vee w^{\circ \circ} \geq \theta\vee w^{\circ \circ}$. Therefore $\varrho\vee w^{\circ \circ} \in \chi$. Then  $\varrho \in E_\chi(W)$. Assume that $\xi, \theta \in E_\chi(W)$. Since $(\xi\wedge \theta) \vee  w^{\circ \circ} = (\xi\vee w^{\circ \circ}) \wedge (\theta\vee w^{\circ \circ}) \in \chi$, then $ \xi \wedge \theta \in\chi$. Clearly, $\chi \subseteq E_\chi(W)$.
\end{proof}
We call $E_\chi(W)$ the extended 
filter of $\chi$. Let $\mathcal{L}$ be any set. Assume that $[0,1]$ is the unit 
interval. A map $\Upsilon:\mathcal{L} 
\rightarrow ([0,1], \wedge, \vee)$ is 
characterised to be a fuzzy subset of $
\mathcal{L}$, with  For every $\theta, \xi \in 
[0,1]$, the operations are defined as follows 
\begin{center}
$\xi \wedge \theta = inf\{\xi, \theta
\}$, \\
$\xi\vee\theta = Sup\{\xi, \theta
\}$.
\end{center}
Throughout the paper we use the notation
$\textbf{Fuzz}(\mathcal{K})$ for the 
collection of every fuzzy subset of a 
bounded distributed lattice $\mathcal{K}$
.
\begin{definition}$[1]$. Let $\Phi,\Psi \in \textbf{Fuzz}(\mathcal{K})$. Define fuzzy subsets $\Phi \cup \Psi$ and $\Phi \cap \Psi$ of $\mathcal{K}$ as follows: let $\xi \in \mathcal{K}$,
\begin{center}
$(\Psi \cup \Phi)(\xi) = \Psi(\xi) \vee \Phi(\xi)$,\\
$(\Psi \cap \Phi)(\xi) = \Psi(\xi) \wedge \Phi(\xi)$.
\end{center}

\end{definition}

\begin{definition}$[2]$
 If $\Upsilon \in \textbf{Fuzz}(\mathcal{K})$. Then $\Upsilon
 $ is characterised as a fuzzy sublattice of $\mathcal{K}$ provided that for every $\theta, \xi \in \mathcal{K}$ the following satisfied;\\ 
(i) $\Upsilon(\xi\wedge \theta)\geq \Upsilon(\xi)\wedge \Upsilon(\theta)$,\\ 
(ii) $\Upsilon(\xi\vee \theta)\geq \Upsilon(\xi)\wedge \Upsilon(\theta)$.\\
\end{definition}

\begin{definition}$[14]$ If $\Upsilon \in \textbf{Fuzz}(\mathcal{K})$, then $\Upsilon$ is a fuzzy ideal if 
the following valid;\\   
(i) $\Upsilon(0)=1$,\\
(ii) For every $\theta, \xi \in \mathcal{K}$, $\Upsilon(\theta\wedge \xi)\geq \Upsilon(\theta)\vee \Upsilon(\xi)$,\\ 
(iii) For every $\theta, \xi \in \mathcal{K}$, $\Upsilon(\theta\vee \xi)\geq \Upsilon(\theta)\wedge \Upsilon(\xi)$.
\end{definition}
We notate the following;
\begin{center}$\textbf{FID}(\mathcal{K}) = \{\Upsilon\in \textbf{Fuzz}(\mathcal{K}), \ 
\Upsilon \ is \ a \ fuzzy  \  ideal \}.$
\end{center} 
\begin{theorem}$[14]$
 If $\Upsilon \in \textbf{Fuzz}(\mathcal{K})$, then $\Upsilon$ 
 is a fuzzy ideal if 
and only if  $\Upsilon(0)=1$ and  
for every $\theta, \xi \in \mathcal{K}$, $\Upsilon(\theta\vee \xi)= \Upsilon(\theta)\wedge \Upsilon(\xi)$.\\
\end{theorem}
\begin{definition}$[14]$ If $\Upsilon \in \textbf{Fuzz}(\mathcal{K})$. Then $\Upsilon$ is a fuzzy filter provided that;\\   
(i) $\Upsilon(1)=1$,\\
(ii) For every $\theta, \xi \in \mathcal{K}$, $\Upsilon(\theta\wedge \xi)\geq \Upsilon(\theta)\wedge \Upsilon(\xi)$,\\ 
(iii) For every $\theta, \xi \in \mathcal{K}$, $\Upsilon(\theta\vee \xi)\geq \Upsilon(\theta)\vee \Upsilon(\xi)$.
\end{definition}
Consider the following set;  
\begin{center}
$\textbf{FFL}(\mathcal{K}) := \{\Upsilon\in \textbf{Fuzz}(\mathcal{K}), \ 
\Upsilon \ is \ a \ fuzzy  \  filter \}.$
\end{center} 
\begin{theorem}$[14]$
 Let $\Upsilon \in \textbf{FFL}(\mathcal{K})$. Then $\Upsilon$ 
 is a fuzzy filter if 
and only if  $\Upsilon(1)=1$ and  
for every $\theta, \xi \in \mathcal{K}$, $\Upsilon(\theta\wedge \xi)= \Upsilon(\theta)\wedge \Upsilon(\xi)$.\\
\end{theorem}
A non constant $\Upsilon \in \textbf{FFL}(\mathcal{K})$ is called proper. A proper 
$\chi \in \textbf{FID}(\mathcal{K})$
is promoted to be a prime fuzzy ideal for
every $\Phi, \Psi \in \textbf{FID}(\mathcal{K})$ satisfying that $\Phi \cap \Psi \subseteq 
\chi$. We have $\Phi \subseteq \chi$ or $\Psi \subseteq \chi$. A proper $\Upsilon \in \textbf{FFL}(\mathcal{K})$ is said to be a prime fuzzy filter if $\Phi, 
\Psi \in \textbf{FFL}(\mathcal{K})$ satisfying that $\Phi \cap \Psi \subseteq 
\Upsilon$. Then $\Phi \subseteq \chi$ or $\Psi \subseteq \Upsilon$.
\section{Extended Fuzzy filters}
For a nonempty subset $W$ of $\mathcal{L}$ and $\chi \in \textbf{FFL}(\mathcal{L})$,  
define 
\begin{center}
    $\Upsilon_{\chi,W}(\theta)= Sup\{ \chi(\theta)\vee \chi(w^{\circ \circ}) ;  \ for \ every \ w \in W \}$,  
\end{center}
for every $\theta \in \mathcal{L}$.
\begin{theorem}
Let $W\subseteq \mathcal{L}$ and $\chi \in \textbf{FFL}(\mathcal{L})$. Then $\Upsilon_{\chi,W}$ is a prime fuzzy filter containing $\chi$.
\end{theorem}

\begin{proof}Clearly, $\Upsilon_{\chi,W}(1)=1$ as $\chi(1)=1$.
We have that \begin{eqnarray*}
  \Upsilon_{\chi,W}(\xi\vee \theta)&=& 
   Sup\{ \chi(\xi \vee \theta)\vee \chi(w^{\circ \circ}) ;  \ for \ 
   every \ w \in W \}\\  &\geq &  Sup\{ (\chi(\xi) \vee 
   \chi(\theta))\vee \chi(w^{\circ \circ}) ;  \ for \ every \ w 
   \in W \}  \\& =&  Sup\{ \chi(\xi )\vee \chi(w^{\circ 
   \circ}) ;  \ for \ every \ w \in W \}  \\&\vee 
   &
   Sup\{ \chi(\theta)\vee \chi(w^{\circ \circ}) ;  \ for \ 
   every 
   \ w \in W \}  \\&= & \Upsilon_{\chi,W}(\xi) \vee 
   \Upsilon_{\chi,W}(\theta).
  \end{eqnarray*} And
  \begin{eqnarray*}
  \Upsilon_{\chi,W}(\xi\wedge \theta)&=& 
   Sup\{ \chi(\xi \wedge \theta)\vee \chi(w^{\circ \circ}) ;  \ for \ 
   every \ w \in W \}\\  &\geq &  Sup\{ 
   (\chi(\xi) \wedge 
   \chi(\theta))\vee \chi(w^{\circ \circ}) ;  \ 
   for \ every \ w \in W \}  \\& =&  Sup\{ 
   \chi(\xi )\vee \chi(w^{\circ \circ}) ;  \ 
   for \ every \ w \in W \}  \\&\wedge 
   &
   Sup\{ \chi(\theta)\vee \chi(w^{\circ \circ}) ;  \ for \ 
   every 
   \ w \in W \}  \\&= & \Upsilon_{\chi,W}(\xi) \wedge 
   \Upsilon_{\chi,W}(\theta).
  \end{eqnarray*}
 Obviously, $\chi \subseteq 
 \Upsilon_{\chi,W}$. Moreover, if $\Phi \cap\Psi 
 \subset \Upsilon_{\chi,W}$. Then $(\Phi \cap 
 \Psi)(\theta) \leq \Upsilon_{\chi,W}(\theta)$. Therefore $
 \Phi(\theta) \wedge \Psi(\theta) \leq \Upsilon_{\chi,W}(\theta)
 $. Hence, $\Phi(\theta)\leq \Upsilon_{\chi,W}(\theta)$ or $ 
 \Psi(\theta) \leq \Upsilon_{\chi,W}(\theta)$.
\end{proof}
We call $\Upsilon_{\chi,W}$ an extended fuzzy filter 
of $\chi$. The following theorem encapsulates some 
basic properties of $\Upsilon_{\chi,W}$. 
\begin{lemma}
For a fuzzy filter $\chi$, the following are valid for a subset $W$ and a member $\theta$ of $\mathcal{L}$
\item{(1)} If $Z \subseteq W$, then $\Upsilon_{\chi,Z} \subseteq \Upsilon_{\chi,W}$,
\item{(2)} If $\chi_1 \subseteq \chi_2$, then $\Upsilon_{\chi_1,W} \subseteq \Upsilon_{\chi_2,W}$,
\item{(3)} The equality $\Upsilon_{\chi,W}(\theta)=\chi(\theta) $ is valid providing that $w^{\circ \circ} \leq \theta$ for every $w \in W$,
\item{(4)} If $\chi$ is a one to one function, then $
\Upsilon_{\chi,W}(\theta)=\chi(\theta) $ implies 
that $w^{\circ \circ} \leq \theta$ for every $
w \in W$,
\item{(5)} If  $ w^{\circ\circ} 
=1 $ for some $ w\in W$, then $ \Upsilon_{\chi,W}
(\theta) = 1$, 
\item{(6)} $\Upsilon_{\chi,\mathcal{L}}(\theta)=1=\Upsilon_{\chi,\{1\}}
(\theta)$. 
\item{(7)} If $\Upsilon_{\chi,W}(\theta)=1$, then either 
$\chi(\theta)=1$ or $\chi(w^{\circ\circ})=1$ for some 
$ w\in W$.
\end{lemma}

\begin{proof}
 (1) If $Z \subseteq W$, 
 we get that $Sup\{ \chi(\theta)\vee 
 \chi(a^{\circ \circ}) ;  \ for \ every \ a \in 
 Z \} \leq 
 Sup\{ \chi(\theta)\vee \chi(w^{\circ \circ}) ;  
 \ for \ every \ w \in W \}$. Hence 
$\Upsilon_{\chi,Z} \subseteq \Upsilon_{\chi,W}$.\\
 (2) Suppose that $\chi_1 \subseteq \chi_2$. We see that $Sup\{ 
 \chi_1(\theta) \vee \chi_1(w^{\circ \circ}) ;  \ for \ every \ w \in 
 W \} \leq Sup\{ \chi_2(\theta)\vee \chi_2(w^{\circ \circ}) ;  \ for 
 \ every \  w \in W \}$. Hence $
 \Upsilon_{\chi_1,W} \subseteq \Upsilon_{\chi_2,W}$.\\
 (3) Let $\theta \leq w^{\circ \circ}$ for 
every $w \in W$. Then $\chi(w^{\circ \circ})  
\leq \chi(\theta)$. Since $\chi(\theta)= \chi(\theta) \vee 
\chi(w^{\circ \circ})$, then 
\begin{eqnarray*}
  \Upsilon_{\chi,W}(\theta)&=& 
   Sup\{ \chi(\theta)\vee \chi(w^{\circ \circ}) ;  \ for \ 
   every \ w \in W \}\\  & = &  Sup\{ \chi(\theta)\}  
   \\& = &   \chi(\theta).
  \end{eqnarray*} 
 (4) Suppose that $\chi(\theta)= \Upsilon_{\chi,W}(\theta)$. 
 Then $\chi(\theta)\geq \chi(w^{\circ \circ})$ for
 every $w \in W$. By the definition of 
 fuzzy filters, for every $w \in W$ we see  
 that $\chi(w^{\circ \circ})= \chi(\theta 
 \wedge w^{\circ  \circ})$. For every $w \in W
 $ and by the one to one property, $
 w^{\circ \circ}=\theta \wedge w^{\circ  
 \circ}$. Hence we get the required result.\\
 (5) If $w^{\circ\circ} =1 $ for some $w
 \in W$, then $\chi(\theta) \vee \chi(w^{\circ 
 \circ})=1$. Hence $\Upsilon_{\chi,W}(\theta)=1$.\\
(6) Follows directly from (5). \\
(7) Let $\Upsilon_{\chi,W}(\theta)=1$. Then $\chi(\theta)\vee 
\chi(w^{\circ \circ})=1$ for some $w\in  
 W$.  As $\chi(\theta), \chi(w^{\circ \circ}) \in [0,
 1]$, hence either $\chi(\theta)=1$ or $\chi(w^{\circ
 \circ})=1$ for some $ w \in W$.
\end{proof}

    \begin{proposition} Let $\mathcal{L} \in 
    \mathbf{MS}$. For nonempty subsets $
    W$ of $
    \mathcal{L}$ and $\chi 
    \in \textbf{FFL}(\mathcal{L})$, we have 
    that
    \item{(1)} $\Upsilon_{\chi_1,W}(\theta)\vee 
    \Upsilon_{\chi_2,W}(\theta) = \Upsilon_{\chi_1\cup \chi_2,W}(\theta)
    $,
     \item{(2)} $\Upsilon_{\chi,W}(\xi)\wedge 
    \Upsilon_{\chi,W}(\theta) = \Upsilon_{\chi,W}
    (\xi\wedge \theta)$.
    \end{proposition}
     \begin{proof}
     \item{(1)} We see that \begin{eqnarray*}
  \Upsilon_{\chi_1,W}(\theta)\vee \Upsilon_{\chi_2,W}(\theta)&=& 
   Sup\{ \chi_1(\theta)\vee \chi_1(w^{\circ \circ}) ;  \ for \ 
   every \ w \in W \}\\  & \vee &  Sup\{ \chi_2(\theta) 
   \vee  \chi_2(w^{\circ \circ});  \ for \ every \ 
   w \in W \}  \\& =&  Sup\{ (\chi_1(\theta )\vee 
   \chi_1(w^{\circ \circ}))\vee (\chi_2(\theta)  \vee  
   \chi_2(w^{\circ \circ})) ;  \ for \ every \ w \in 
   W \}  \\&= &Sup\{ (\chi_1(\theta )\vee \chi_2(\theta)) 
   \vee (\chi_1(w^{\circ \circ}) \vee  \chi_2(w^{\circ 
   \circ})) ;  \ for \ every \ w \in 
   W \}  \\&= &Sup\{ ((\chi_1\vee\chi_2)(\theta)) 
   \vee ((\chi_1 \vee \chi_2)(w^{\circ 
   \circ})) ;  \ for \ every \ w \in 
   W \}  \\&= & \Upsilon_{\chi_1\cup \chi_2,W}(\theta).\\
   \end{eqnarray*}
    \item{(2)} By Theorem $(3.1)$, we get the required. 
    
     \end{proof}
 
\begin{definition}
Let $\chi\in \textbf{FFL}(\mathcal{L})$. We say 
that $\chi$ is a fixed fuzzy filter relative to a
non empty set $W$ if for every $
\theta$ the equality $\Upsilon_{\chi,W}(\theta)=\chi(\theta)$ is satisfied. 
 \end{definition}


  \begin{example}
  \item{(1)} In the virtue of Lemma $3.2(3)$, 
  every  $\chi\in \textbf{FFL}(\mathcal{L})$ is 
  a fixed fuzzy filter relative to $W=\{0\}$. 
  That is $\Upsilon_{\chi,\{0\}}(\theta)=\chi(\theta)$.
  \item{(2)} Let $A=\{\rho\in \mathcal{L} \ ; 
  \rho^{\circ \circ}=0\}$. Then $
  \Upsilon_{\chi,A}(\theta)=\chi(\theta)$. Since $\rho^{\circ 
  \circ}=0\leq \theta$ for every $\theta \in \mathcal{L}$. Thus 
  $\chi(\rho^{\circ  \circ})\leq \chi(\theta)$. 
  \item{(3)} Let $C=\{\rho\in \mathcal{L} 
  \ ; \chi(\rho^{\circ \circ})=0\}$. By 
  the definition
 of $\Upsilon_{\chi,C}$, we get for every $\theta \in  \mathcal{L} $ that $ \Upsilon_{\chi,C}(\theta)=\chi(\theta)$. 
  \item{(4)} Consider the following Hasse 
  diagram 
 \newline 
 \begin{center}
   \begin{tikzpicture}[scale=.5]
  \node (o) at (0,2) {$1$};
\node (a) at (2,0) {$\theta =\theta^{\circ}$};
  \node (b) at (-2,0) {$\xi^{\circ}=\xi$};
  \node (z) at (0,-2) {$0$};
\draw (o) -- (a) -- (z) -- (b) -- (o)  ; 
\end{tikzpicture}
 \end{center} 
 \begin{center}
    Figure(1): $\mathcal{L}$
\end{center} 
Obviously, $(\mathcal{L},^\circ)\in\mathbf{MS}
$. Define a fuzzy subset $\chi$ of $\mathcal{L}
$ as $\chi(\theta)=\chi(1)=1$ and $\chi(0)=
\chi(\xi)=0.5$. 
Then $\chi\in \textbf{FFL}(\mathcal{L})$. Let 
$W=\{0,\xi \}$. Obviously, $\chi$ is a fixed 
fuzzy filter relative to the set $W$.\\ 
\end{example}

\begin{proposition}
Let $Z,W \subseteq \mathcal{L}$. If a 
fuzzy filter $\chi$ is 
fixed relative to $W$ and $Z\subseteq W$, 
then $\chi$ is a fixed fuzzy filter relative to $Z$.
\end{proposition}
\begin{proof}
 Let $Z \subseteq W$. Then $\Upsilon_{\chi,Z}(\xi) \subseteq \Upsilon_{\chi,W}(\xi)=\chi(\xi)$. Therefore $ \Upsilon_{\chi,Z}(\theta)=\chi(\theta)$. Hence $\chi$ is a fixed fuzzy filter relative to $Z$.
\end{proof}
\begin{proposition}
If $\chi_1,\chi_2$ are fixed fuzzy filters 
relative to a subset $W\subseteq \mathcal{L}$, 
then $\chi_1\cup \chi_2$ is a fixed fuzzy 
filters relative to $W$.
  
\end{proposition}
  \begin{proof}
  This follows easily from Proposition $(3.3)$.  
  \end{proof} 
   \begin{theorem}
     The image of an element $\theta\in \mathcal{L}$ under the 
     fuzzy filter $\Upsilon_{\chi,\{w\}}$ is $\chi(\theta)
     $ or $\chi(w^{\circ \circ})$ for every $w\in W$.
     \end{theorem}
     \begin{proof}One can see that
     \begin{eqnarray*}
  \Upsilon_{\chi,\{w\}}(\theta)&=& 
   Sup\{ \chi(\theta)\vee \chi(w^{\circ \circ}) ;  \ for \ 
   every \ w \in W \}\\  & =&  Sup\{ (\chi(\theta )\vee 
   \chi(w^{\circ \circ})\}  \\&= & \chi(\theta)
   \vee \chi(w^{\circ \circ}) .
   \end{eqnarray*}
     \end{proof}
     This gives immediately the following
     results.
     \begin{corollary}
     For $w \in \mathcal{L}$, if $\Upsilon_{\chi,\{w\}}(\theta)
     \not = \chi(w^{\circ \circ})$, then $\chi$ is a 
     fixed fuzzy filter with respect to $\{w\} 
     $.
     \end{corollary}
     \begin{corollary}
     If $\chi\in \textbf{FFL}(\mathcal{L})$ and 
     $w \in \mathcal{L}$, then $\Upsilon_{\chi,\{w\}}(w)  
     =\chi(w^{\circ \circ})$.
     \end{corollary}
    \section{Strong Fuzzy Extensions} 
\hspace{.5cm}In this section, a more generalized case of extended 
 fuzzy filter called strong fuzzy extension is studied.
 \begin{definition}
 Let $\chi\in \textbf{FFL}(\mathcal{L})$. Define the 
 following set
 \begin{center}
 $\Omega_{\chi,W}(\theta)=Sup \{\chi(\theta\vee w^{\circ \circ}) 
 ;  \ for \ every \ w \in W  \}$.
 \end{center}
 \end{definition}
 We see that $\Omega_{\chi,W}$ is an 
 extension of $\Upsilon_{\chi,W}(\theta)$. By using the 
 properties of fuzzy filters, we see that 
for every $\xi, \theta \in \mathcal{L}$, $\chi(\xi\vee \theta)\geq \chi(\xi)
\vee \chi(\theta)$. Therefore, $\Omega_{\chi,W}(\theta)\geq 
\Upsilon_{\chi,W}(\theta)$. Hence $\Upsilon_{\chi,W}
\subseteq\Omega_{\chi,W}$. The inclusion $\Upsilon_{\chi,W}
\subseteq \Omega_{\chi,W}$ is proper as verified in the next example.
     \begin{example}
 Consider the lattice with the following Hasse diagram $\mathcal{L}$ in Figure (1).
 \newline 
 \begin{center}
   \begin{tikzpicture}[scale=.5]
  \node (max) at (0,8) {$1$};
  \node (u) at (0,5) {$u$};
\node (z) at (0,2) {$z$};
\node (y) at (2,0) {$\xi$};
  \node (x) at (-2,0) {$\theta$};
  \node (t) at (0,-2) {$t$};
  \node (mn) at (0,-4) {$0$};

\draw (mn) -- (t) -- (y) -- (z) -- (u) -- (max) ;
\draw (t) -- (x) -- (z); 
\end{tikzpicture}
 \end{center}
\begin{center}
    Figure (1): $\mathcal{L}$
\end{center} 

 Define a unary operation $^{\circ}$ on $
 \mathcal{L}$ by 
$\theta^{\circ}=t,\xi^{\circ}=z^{\circ}=t^{\circ}=u,u^{\circ}=\xi
, 1^{\circ}=0,0^{\circ}=1$. Then $
(\mathcal{L},^\circ)\in \mathbf{MS}$. Define 
the following fuzzy set $\chi:\mathcal{L}\longrightarrow [0,1]$ 
as the following; $\chi(0)=0, \ \chi(u)=0.8, \ 
\chi(t)=\chi(\xi)=\chi(\theta)=0.6, \ \chi(z)=0.7$ and $\chi(1)=0.7$. 
Obviously, $\chi \in \textbf{FFL}(\mathcal{L})$. Then
\begin{eqnarray*}
    \Upsilon_{\chi,\{\xi\}}(\theta)&=& Sup\{\chi(\theta)\vee 
    \chi(\xi^{\circ \circ})\}\\  &=& \chi(\theta)\vee 
    \chi(\xi^{\circ \circ})\\&=& \chi(\theta)\vee 
    \chi(\xi)=0.6.
   \end{eqnarray*}
   And
    \begin{eqnarray*}
     \Omega_{\chi,\{\xi\}}(\theta)&=& Sup\{\chi(\theta\vee 
    \xi^{\circ \circ})\}\\  &=& \chi(\theta\vee 
    \xi^{\circ \circ})\\&=& \chi(\theta\vee 
    \xi)=0.7.
    \end{eqnarray*}
  \end{example}
   We call $\Omega_{\chi,W}$ a strong 
   fuzzy extension of $\chi$.
    \begin{theorem}
Let $\chi \in \textbf{FFL}(\mathcal{L})$. If $W \subseteq \mathcal{L}
$, then $\Omega_{\chi,W}(\theta)\in \textbf{FFL}(\mathcal{L})$. 
\end{theorem}

\begin{proof}
Obviously, $\Omega_{\chi,W}(1)=1$. 
We have that \begin{eqnarray*}
  \Omega_{\chi,W}(\xi\vee \theta)&=& 
   Sup\{ \chi((\xi \vee \theta)\vee w^{\circ \circ})) ;  \ for \ 
   every \ w \in W \}\\  &= &  Sup\{ \chi((\xi \vee 
   w^{\circ \circ}) \vee (\theta\vee w^{\circ 
   \circ}))  ;  \ for \ every \ w 
   \in W \}  \\& \geq&  Sup\{ \chi(\xi\vee w^{\circ 
   \circ}) \vee \chi(\theta\vee w^{\circ 
   \circ}) ;  \ for \ every \ w \in W \}  \\&= 
   &Sup\{ \chi(\xi\vee w^{\circ \circ}) ;  \ for 
   every  \ w \in W \}  \\&\vee &
   Sup\{ \chi(\theta\vee w^{\circ \circ}) ;  \ for \ 
   every \ w \in W \}  \\&= & \Omega_{\chi,W}(\xi) 
   \vee 
   \Omega_{\chi,W}(\theta).
  \end{eqnarray*}
  \begin{eqnarray*}
  \Omega_{\chi,W}(\xi\wedge \theta)&=& 
   Sup\{ \chi((\xi \wedge \theta)\vee w^{\circ \circ})) ;  \ for \ 
   every \ w \in W \}\\  &= &  Sup\{\chi((\xi \vee 
   w^{\circ \circ}) \wedge (\theta\vee w^{\circ 
   \circ}))  ;  \ for \ every \ w 
   \in W \}  \\& \geq&  Sup\{ \chi(\xi\vee w^{\circ 
   \circ}) \wedge \chi(\theta\vee w^{\circ 
   \circ}) ;  \ for \ every \ w \in W \}  \\&= 
   &Sup\{ \chi(\xi\vee w^{\circ \circ}) ;  \ for 
   every  \ w \in W \}  \\&\wedge &
   Sup\{ \chi(\theta\vee w^{\circ \circ}) ;  \ for \ 
   every \ w \in W \}  \\&= & \Omega_{\chi,W}(\xi) 
   \vee 
   \Omega_{\chi,W}(\theta).
  \end{eqnarray*}
\end{proof}
\begin{remark}
   We see that $\Upsilon_{\chi,W}$ 
   is not necessarily equal to $\Omega_{\chi,
   W}$. But the 
   equality occurs if $\chi$ satisfying that 
   for every $\xi, \theta \in \mathcal{L}$, $
   \chi(\xi \vee \theta)=\chi(\xi)\vee 
   \chi(\theta)$. i.e $\Omega_{\chi,W}
   (\theta) =\Upsilon_{\chi,W}
   (\theta)$ 
   for every $\theta \in \mathcal{L}$.
\end{remark}
    For a fuzzy set $\mu$ in $S$ and 
    $t \in [0,1]$, we set  $\mu_t:=\{\theta\in S;\ \mu
    (\theta)\geq t\}$.
    
    \begin{definition}
    An element $\alpha$ is said to be a dense 
    element in a set $W$ with 
    respect to fuzzy  
    filter $\mu$ if $\alpha \in \mu_t$ such  
    that $t=\max \{\mu(\theta); \theta\in 
    W\}$.  
    \end{definition}
    \begin{example}
    Let $\mathcal{L}=\aleph$ (Chain of natural numbers). 
    Define a fuzzy set $\lambda:\mathcal{L}
    \longrightarrow [0,1]$ such that $
    \lambda(n)=\dfrac{1}{n}$. Obviously, $1$ is 
    a dense element with respect to any subset 
    of $\aleph$.  
    \end{example}
    The symbol $W^{\circ \circ}$ stands 
    for the set $\{w^{\circ \circ}; w\in 
    W\}$. The following theorem 
    illustrates the necessity of dense 
    element. 
    \begin{theorem}
    If $w^{\circ \circ}_\theta$ is a dense 
    element in the set $W^{\circ \circ}$
    with respect to a fuzzy filter $\chi$. 
    Then $\Upsilon_{\chi,W}(\theta)= 
    \chi(\theta)\vee \chi(w^{\circ \circ}_\theta)$.
    \end{theorem}
    \begin{proof}
    By the definition of dense elements,  
    we have that $\chi(w^{\circ \circ}
    _\theta) \geq \chi(w^{\circ \circ})$
    for every $w^{\circ \circ}\in W^{\circ 
    \circ}$. Hence $\Upsilon_{\chi,W}
    (\theta)= \chi(\theta)\vee 
    \chi(w^{\circ \circ}_\theta)$.
    \end{proof}
      The following theorem characterises 
    the value of the strong fuzzy extended 
    filter in terms of a dense element.  
   
      \begin{theorem}
    Let $w\in W$, the equality $
    \Omega_{\chi,W}(\theta)
    =\chi(\theta\vee w^{\circ \circ})$ is 
    valid if and only if $
    \theta\vee w^{\circ \circ}$ is a dense 
    element in the set $ \{\theta\vee 
    w^{\circ \circ};  w\in W\}$ with respect to a fuzzy filter $\chi$. 
    \end{theorem}
    \begin{proof}
    Consider $\Omega_{\chi,W}
    (\theta)=\chi(\theta\vee w^{\circ 
    \circ})$, thus the element $\theta\vee 
    w^{\circ \circ}$ has the largest value 
   among all such members of the set $
   \{\theta\vee  w^{\circ \circ};  w\in 
   W\} $ under the fuzzy set $\chi$. The 
   other direction follows directly.
    \end{proof}
    \section{Homomorphisms of fuzzy extensions}
   \begin{theorem}
   Let $\chi \in \textbf{FFL}(\mathcal{L})$ and $W\in \mathcal{L}$.
   If $\chi$ is a join homomorphism map, then $
   \Upsilon_{\chi,W}$ is a lattice homomorphism.
   Furthermore, if $\chi$ is an $\emph{MS}
   $-homomorphism, then so $\Upsilon_{\chi,W}$.  
   \end{theorem}
   \begin{proof}
   Every fuzzy filter $\mu$ have the 
   property that $\mu(\theta\wedge \xi)=\mu(\theta)\wedge \mu(\xi)$. 
   Consequently, $\Upsilon_{\chi,W}(\theta\wedge \xi)=
   \Upsilon_{\chi,W}(\theta)\wedge \Upsilon_{\chi,W}(\xi)$. 
   Assume that $\chi$ is a join homomorphism. Then 
   $\chi(\xi\vee\theta)=\chi(\xi)\vee \chi(\theta)$. Therefore
   \begin{eqnarray*}
  \Upsilon_{\chi,W}(\xi\vee \theta)&=& 
   Sup\{\chi((\xi\vee \theta)\vee w^{\circ \circ})) ;  
   \ for \  every \ w \in W \} \\  &= &  Sup\{\chi(\xi)\vee \chi(\theta) \vee \chi(w^{\circ \circ})  ;  \ 
   for \ every \ w \in W \}  \\&=&  Sup \{(\chi(\xi)
   \vee \chi(w^{\circ \circ})) \vee (\chi(\theta)\vee 
   \chi(w^{\circ \circ})) ;  \ for \ every \ w \in 
   W \}  \\&= 
   &Sup\{\chi(\xi)\vee \chi(w^{\circ \circ}) ;  \ for 
   \ every  \ w \in W \}  \\&\vee &
   Sup\{\chi(\theta)\vee \chi(w^{\circ \circ}) ;  \ for \ 
   every \ w \in W \}  \\ &= & \Upsilon_{\chi,W}(\xi) 
   \vee  \Upsilon_{\chi,W}(\theta).  \end{eqnarray*}
  Furthermore, for every $\theta\in \mathcal{L}$ suppose that  
  $\chi(\theta^{\circ \circ})= (\chi(\theta))^{\circ \circ}$. 
  Then;
  \begin{eqnarray*}
  \Upsilon_{\chi,W}(\theta^{\circ \circ})&=& 
   Sup\{(\chi(\theta^{\circ \circ})\vee \chi(w^{\circ 
   \circ}) ;   \ for \  every \ 
   w \in W \}\\  &= &  Sup\{(\chi(\theta))^{\circ 
   \circ}\vee (\chi(w^{\circ 
   \circ}))^{\circ \circ} ; \ for \ every \ w 
   \in W \}  \\ &= &  Sup\{(\chi(\theta)\vee \chi(w^{\circ 
   \circ}))^{\circ \circ} ; \ for \ every \ w 
   \in W \}  \\ &=&  (\Upsilon_{\chi,W}(\theta))^{\circ 
   \circ}.   
  \end{eqnarray*}
  
   \end{proof}
   \begin{proposition}
   Let $\chi\in \textbf{FFL}(\mathcal{L})$ 
   and $W
   \subseteq \mathcal{L}$. Every element 
   in the kernel of
   $\Upsilon_{\chi,W}$ is an element of 
   the kernel of $\chi$ and $W^{\circ 
   \circ}$ is a
   subset of the kernel of $\chi$. 
   Furthermore, the converse is true.
   \end{proposition}
   \begin{proof}
   Assume that $\theta$ is an element of 
   the kernel of 
   $\Upsilon_{\chi,W}(\theta)$. 
   Then $\chi(\theta)\vee 
   \chi(w^{\circ \circ})=0$ for every $w 
   \in W$. 
   Thus $\chi(\theta)=0$ and $
   \chi(w^{\circ \circ})=0$ for 
   every $w \in W$. Hence $\theta$ and 
   $w^{\circ \circ}$ are elements of the 
   kernel of $\chi$ for every $w \in W$. 
   Conversely, let $\theta$ be an element 
   of the kernel of $\chi$ and let
   $W^{\circ \circ}$ be a subset of the 
   kernel of $\chi$. Then $\chi(\theta)=0$ 
   and $\chi(w^{\circ \circ})=0$ for every 
   $w \in W$. Hence $\theta$ be 
   an element of the kernel of $
   \Upsilon_{\chi,W}$. 
   \end{proof}
   \begin{proposition}
   An element is in the cokernel of
   $\Upsilon_{\chi,W}$ if and only if it 
   is an element of the cokernel of $\chi$ 
   or the intersection between $W^{\circ 
   \circ}$ and the cokernel of $\chi$ is 
   not empty. 
   \end{proposition}
   \begin{proof}
   Suppose that $\theta$ is an element in 
   the cokernel of $\Upsilon_{\chi,W}$. 
   Then $ \Upsilon_{\chi,W}(\theta)=1$. 
   Thus $\chi(\theta)\vee 
    \chi(w^{\circ \circ})=1$ for some $w
    \in W$.  As $\chi(\theta), \chi(w^{\circ \circ}) 
    \in [0,1]$, hence either $\chi(\theta)=1$ or 
    $\chi(w^{\circ\circ})=1$ for some $ w 
    \in W$. The reverse implication follows directly 
    from the definition of $\Upsilon_{\chi,W}$.
   \end{proof}
  Let $\theta \in \mathcal{L}$. Define the following set
   \begin{center}
 $\underleftarrow{\Upsilon}(\Upsilon_{\chi,W}(\theta))=
\lbrace a\in \mathcal{L}; \ \Upsilon_{\chi,W}(a)=
\Upsilon_{\chi,W}(\theta) \rbrace$.
   \end{center}
 The set $\underleftarrow{\Upsilon}
 (\Upsilon_{\chi,W}(\theta))$ is the inverse (as a 
 relation) of $\Upsilon_{\chi,W}(\theta)$.
  \begin{lemma}
  For every $\theta \in \mathcal{L}$. The set $
  \underleftarrow{\Upsilon}(\Upsilon_{\chi,W}(\theta))$ 
  is closed under the meet operation. 
  Furthermore, If $\chi$ is a join homomorphism, then $
  \underleftarrow{\Upsilon}(\Upsilon_{\chi,W}(\theta))$
  is a subalgebra of $\mathcal{L}$. 
  \end{lemma}
  \begin{proof}
  Assume that $\xi,\rho \in 
  \underleftarrow{\Upsilon}(\Upsilon_{\chi,W}(\theta))
  $. Therefore
   \begin{eqnarray*}
  \Upsilon_{\chi,W}(\xi\wedge \rho)&=& 
   \Upsilon_{\chi,W}(\xi) \wedge \Upsilon_{\chi,W}
   (\rho)\\  
   &= & \Upsilon_{\chi,W}(\theta) \wedge \Upsilon_{\chi,W}
   (\theta)  \\&= &  \Upsilon_{\chi,W}(\theta).
  \end{eqnarray*}
  Hence $\xi\wedge \rho \in \Upsilon_{\chi,W}(\theta)
  $. Let $\chi$ be a $\vee$-homomorphism 
  . By Theorem $(5.1)$, $\Upsilon_{\chi,W}
  $ is a lattice homomorphism. Let $\xi,\rho \in \underleftarrow{\Upsilon}(\Upsilon_{\chi,W}(\theta))
  $. Then 
  \begin{eqnarray*}
  \Upsilon_{\chi,W}(\xi\vee \rho)&=& 
   \Upsilon_{\chi,W}(\xi) \vee \Upsilon_{\chi,W}
   (\rho)\\  
   &= & \Upsilon_{\chi,W}(\theta) \vee \Upsilon_{\chi,W}
   (\theta)  \\&= &  \Upsilon_{\chi,W}(\theta).
  \end{eqnarray*}
  
  \end{proof}
  
  \section*{Statements and Declarations}
 \subparagraph{\textbf{Confict of interest} 
} 
 The authors declare that they have no conflict of interest.\\
 \subparagraph{\textbf{Competing Interests} 
} 
This research did not receive any specific grant from funding agencies in the public, commercial, or not-for-profit sectors.\\
 \subparagraph{\textbf{Author Contributions} 
} 
 All authors contributed to the study conception of this paper.

\end{document}